\documentclass[12pt]{amsart}
\newtheorem{thm}{Theorem}
\newtheorem{lem}[thm]{Lemma}
\newtheorem{de}[thm]{Definition}
\newtheorem{prop}[thm]{Proposition}

\def\N{{\Bbb N^+}}

\def\nt{\noindent}
\def\piff{\leftrightarrow}

\def\lf{\lfloor}
\def\rf{\rfloor}
\def\ph{\varphi}

\def\v{\Vert}
\def\be{\begin{enumerate}}
\def\ee{\end{enumerate}}
\def\ar{\rightarrow}
\def\to{\rightarrow}

\def\Cal{\mathcal}

\def\amko{\alpha(\langle m,k\rangle)}
\def\amk{\exists m\,\amko}
\def\apko{\alpha(\langle p,k\rangle)}
\def\apk{\exists p\,\apko}

\def\amnk{\exists m\leq n\,\amko}
\def\bnk{\exists n\,\beta(n)=k}

\def\vv{\,\vee\,}
\def\pnp{\ph\vee\neg\ph}

\def\er{\exists^R}
\def\fr{\forall^R}
\def\en{\exists^{\N}}
\def\fn{\forall^{\N}}

\begin{document}
\vskip5cm
\nt {\large\bf Generalized Kripke's Schema and the Expressive Power of Intuitionistic Real Algebra}

\vskip1cm

\nt{\bf Mikl\'os Erd\'elyi-Szab\'o}
\vskip0.5cm

\nt R\'enyi Alfr\'ed Institute of Mathematics,

\nt 1053 Budapest, Re\'altanoda u. 13-15.

\nt Hungary\footnote{e-mail: mszabo@@renyi.hu}
\vskip1cm

\nt {\bf Abstract.}\quad We introduce a relativized version of random Kripke's schema  and show how it may be applied in the investigation of the expressive power of intuitionistic real algebra by interpreting second-order Heyting arithmetic in it. 
\vskip0.6cm

\nt {\bf Mathematics Subject Classification:} 03-B20, 03-F25, 03-F35, 03-F60.
\vskip0.6cm

\nt {\bf Keywords:} Intuitionism, Kripke's Schema, Second-order Heyting arithmetic, Real algebra, Interpretation.
\pagestyle{myheadings}
\markboth{Mikl\'os Erd\'elyi-Szab\'o}{Generalized Kripke's schema} 
\vskip1cm

\section{Introduction}
To conform with the framework used in previous papers \cite{sco1}, \cite{sco2},\cite{eszmjsl} and \cite{eszmmlq2}, we shall use the language $L_1$ from \cite{sco1}.
It contains 
two sorts of variables -- $m$, $n$, $k$, etc. ranging over the elements of 
$\omega$, and $\alpha$, $\beta$, etc. ranging over choice sequences.
The language contains symbols for certain primitive recursive 
functions like the pairing function $\langle.,.\rangle$ and relations defined on the elements of $\omega$ 
and we also have the equality symbol $=$.
 It 
will be used in atomic formulas of the form 
$t=t'$ or $\xi(t)=t'$ where $t$ and $t'$ are terms of natural-number sort and 
$\xi$ is of choice sequence sort.

\smallskip

In \cite{eszmmlq2} we defined randomized Kripke's schema as the following axiom schema of second-order arithmetic:

\[ \hbox{R-KS}(\ph)\equiv \exists\beta[(\exists n(\beta(n)>0)\ar\ph)\wedge 
(\neg\exists n(\beta(n)>0)\ar\neg\ph)\wedge \]
\[\forall k>0(\neg\exists n(\beta(n)=k)\ar\ph\vee\neg\ph)\wedge 
\forall k>0\,\forall n((\beta(n)=k)\ar\forall m\geq n(\beta(m)=k))]\]

\nt where $\ph$ is a formula that does not contain the choice sequence variable $\beta$ free. 

\medskip

\nt Note that from R-KS the Kripke's schema 

\[\hbox{KS}(\ph)\equiv \exists\beta[\exists n(\beta(n)>0)\ar\ph)\wedge 
(\neg\exists n(\beta(n)>0)\ar\neg\ph\] obviously follows.

\medskip

\nt Also in \cite{eszmmlq2} we have proved the following.
\begin{thm}
\begin{enumerate}
\item[(i)] The models of intuitionistic second order arithmetic described in \cite{sco1} and in \cite{sco2} are models of R-KS. 
\item[(ii)] From a standard axiom system augmented with R-KS the definability in the language of ordered rings of the set of natural numbers follows. 
\end{enumerate}
\end{thm}

\nt In this paper we extend this result by giving a relativised version R-KS (called RR-KS) below
and as an application of this new schema an interpretation of second-order Heyting arithmetic ${\bf HAS}$  in intuitionistic real algebra.

\smallskip

A faithful translation from the usual language of ${\bf HAS}$ containing species (set of natural numbers) variables but no sequence variables into $L_1$ is given in \cite{vandalen2} as mentioned also in \cite{sco1}. If $\alpha$ is a choice sequence, $\{x\,|\,\exists y\, \alpha(\langle x,y\rangle)\neq 0)\}$ is a species of natural numbers and Kripke's schema together with relativised dependent choice $RDC$ 
%(true in the models we shall deal with in Section \ref{models} below) 
imply that any definable species of natural numbers has this form. In what follows we shall use the two versions (using species sort or sequence sort) interchangeably.

\section{RR-KS - Relativised Random Kripke's Schema and its justification using a creative subject argument}

\[\hbox{RR-KS}(\ph,\alpha)\equiv\exists\beta[((\exists n\,\beta(n)>0\ar\ph)\,\wedge\]
\[(\forall k>0\,\amk=1\,\wedge\, \neg\exists n\,\beta(n)>0) \,\ar \neg\ph\,\wedge\]
\[\forall k>0 \,(\neg\bnk\to\pnp)\to (\amk=1\vee\pnp)\,\wedge\]
\[\forall k>0 \,(\amk=1\to (\neg\bnk\to\pnp)\,\wedge\]
\[\forall k>0\,\forall n\,(\beta(n)=k\ar\forall m\geq n\, \beta(m)=k))]\]
where $\ph$ is a formula that does not contain the choice sequence variable $\beta$ free, $\alpha$ is a choice sequence parameter with values in $\{0,1\}$. 
 \medskip
Note that R-KS$(\ph)$ is equivalent to RR-KS$(\ph,\alpha)$ when $\forall k>0\, \amk = 1$. 

\medskip
%\subsection
\noindent {\bf An argument in favor of RR-KS}\label{just}
is a modified version of our argument for R-KS in \cite{eszmjsl}.

\nt Let us suppose that the formula $\ph$ and the choice sequence $\alpha$ are given. 
\nt We fix some sequence of moments of time $0$, $1$, $2$, etc.
 At moment $n$ 

\be
\item for each $p,k\leq n$ we check if $\apko = 1$, 
\item  we also check if we have a proof of $\ph$ or a proof of $\neg\ph$. 
\ee

\nt The values $\beta(n)$ are 
obtained in the following way.

\nt If at moment $n$ 
$\ph$ has not been proved yet, we set $\beta(n)=0$.
Otherwise we chose a random number $0< k\leq n$, and check whether $\amnk=1$. If this is the case, we set  $\beta(n)=k$, if not (ie. $\ph$ is not proved, or $\neg\amnk = 1$), we set $\beta(n)=0$. 
The point here is that there is no information about $k$ 
before $\ph$ has been proved. 
If $\beta(n)=k$ for  some $n$ and $k>0$, then let $\beta(m)=k$ for every $m\geq n$ making
 the last conjunct in the matrix of the schema $RR-KS$ true. 
If at some moment $\neg\ph$ is proved, we set $\beta(n)=0$ for every $n$. 
As usual, we assume that at each moment $n$ either $\ph$ is proven, or it is not, 
and that we have a method to check the value of $\amko$, 
so we can use the above case analysis, and 
$\beta(n)=0\vee\neg\beta(n)=0$ is true for every $n$. 

The first conjunct $\exists n\,\beta(n)>0\ar\ph)$ is obvious, until $\ph$ is proved, the value of $\beta$ is $0$. 

For the second conjunct note, that if $\forall k>0\,\amk=1$, then $\neg\exists n\,\beta(n)>0$ if and only if $\ph$ is never proved. If $\ph$ is true, then it is plausible to assume that it is impossible that $\ph$ will never be proved, so $\ph\ar\neg\neg\exists n\,\beta(n)>0$ which is equivalent to $\neg\exists n\,\beta(n)>0\ar\neg\ph$. 

\nt Next we argue that $\forall k>0 \,(\neg\bnk\to\pnp)\to (\apk=1\vee\pnp)$.

\smallskip 

Assume that for some $k$ $\neg\bnk\to\pnp$ is known. We have the following possibilities.
\be
\item $\bnk$ is known. By the construction of $\beta$ this can happen only if $\ph$ is proven.

\item $\neg\bnk$ is known. Then by the assumption we have $\ph\vee\neg\ph$.
\item None of the above, ie. we have no information on the future values of $\beta$, though we do know that $\neg\bnk\to\pnp$. We have the following subcases.
\be
\item $\apk = 1$ is known. Then we are done. Note also, that if in this case $\neg\bnk$ turns out to be true, this can happen only when $\ph$ or $\neg\ph$ is known.
\item If $\neg\apk = 1$ were known, we would have $\neg\bnk$. 
\item If neither $\apk = 1$, nor $\neg\apk = 1$ is known, and we have no information on  $\ph$ or $\neg\ph$, then we have no information on the values of  $\beta$. In particular  in this case $\neg\bnk\to\pnp$ cannot be known.
\ee
\ee

\smallskip

The next conjunct to check is
$\forall k>0 \,(\amk=1\to (\neg\bnk\to\pnp)$, 

\nt so assume $\amk=1$ for some $k>0$. Then we can find the least such $m$ and at that moment $\exists p\leq m\,\alpha(\langle p,k\rangle)=1)$. Now assume $\neg\bnk$. By the way $\beta$ was constructed,  this can happen only in two ways. Either $\ph$ is proved and the random number picked is not $k$, or $\neg\ph$ is proved. Either way  $\neg\bnk\to\pnp$ follows.

\smallskip
The last conjunct is immediate.

\section{An application: interpretation of intuitionistic second order arithmetic in intuitionistic real algebra}
\vskip0.8cm

\nt Let $L$ denote the language of ordered rings. In \cite{me}, using a coding  similar to the one used by {\sc Cherlin} in \cite{cherlin}, we showed that the $L$-theory of Scott's topological model for intuitionistic analysis is undecidable by encoding true first-order arithmetic in that structure. 
({\sc Scott} defined the model in question in \cite{scott}.) 
In \cite{eszmjsl} we showed that true (classical) first-order arithmetic is interpretable in a class of models which includes the well-known topological models as well as Scowcroft's model (defined in \cite{sco1}) and its generalizations, in \cite{eszmmlq2} we extended this result by showing that classical second order arithmetic is interpretable in intuitionistic real algebraic structures. Here we show {\it intuitionistic} second order arithmetic is interpretable in intuitionistic real algebra in the presence of RR-KS. 
For the sake of convenience, we restrict ourselves to the set of non-negative reals, so in what follows, ``real'' will mean ``non-negative real''. 
\smallskip

First, let us see how the algebraic structure of reals is defined in intuitionistic second order arithmetic.  
\medskip

In~\cite[pages 134-135]{kv} Vesley considers a species $R$ of {\it real-number generators}: $\xi\in R$ (also denoted by $R(\xi)$)
if and only if the sequence $2^{-x}\xi(x)$ ($x\in\omega$) of diadic fractions is a Cauchy-sequence with
$\forall k\exists x\forall p |2^{-x}\xi(x)-2^{-x-p}\xi(x+p)|<2^{-k}$, i.e.
if and only if $\forall k\exists x\forall p 2^k|2^p\xi(x)-\xi(x+p)|<2^{x+p}$.

Equality, ordering, addition and multiplication on $R$ are defined as follows.
\be
\item $\xi=\eta$ if and only if $\forall k\exists x\forall p 2^k|\xi(x+p)-\eta(x+p)|<2^{x+p}$, 
\item $\xi<\eta$ if and only if 
$\exists k\exists x\forall p 2^k(\eta(x+p)\dot -\xi(x+p))\geq 2^{x+p}$, 
\item $(\xi+\eta)(x):=\xi(x)+\eta(x)$ and 
\item $(\xi\eta)(x):=\lf 2^{-x}\xi(x) \eta(x)\rf$.
\ee 
The following facts are also proved in~\cite{kv}. 
If $\xi$, $\eta\in R$ then $\xi+\eta\in R$, $=$ is a congruence relation with respect to $<$ and $+$. 
Similar facts are true for multiplication as well (cf. also \cite[pages 20-21]{heyting}).
\medskip

Using the facts mentioned above, these definitions can be extended readily to polynomials of choice sequences.

$\xi$ is a {\it global real-number generator} just in case $ R(\xi)$ holds. We shall use the letters  $f$, $g$, $u$ etc.  to range over global real-number generators and we shall use the defined quantifiers 
$$\er u\,\theta:\equiv  \exists u(R(u)\wedge \theta)$$
$$\fr u\,\theta:\equiv  \forall u(R(u)\ar \theta)$$

For each natural number $n$ there is a corresponding global real-number generator $f_n$ defined as follows: 
$f_n(l) = m$ iff $m=n2^l$. Then 
\[  f_n f=\overbrace{f+\cdots+f}^{n} \]
If it does not cause any confusion, we shall denote this particular global real-number generator by $n$. 
This way we defined a model for the language of ordered rings with addition and multiplication defined above and the interpretation of 
$0$ and $1$ is $f_0$ and $f_1$ respectively. Note that $=$ has the usual properties, e.g.  $f=g\piff \neg(g<f\lor f<g)$,  $f=g\wedge\ph(f)\ar\ph(g)$ 
etc. $f\not=g$ is defined as $f<g\vee g<f$. Finally note that in general $\neg f=g\not\piff f\not=g$.

\vskip0.8cm
 
In \cite{eszmmlq2} we defined the notion of NE-quotient (non-excluded quotient) of two real number generators $u$, $v\neq 0$ of a model with respect to an element $B$ of the underlying Heyting algebra as a positive number $n$ with the property that $\v\neg n \, v = u\v\leq B$. Here we need a modified version (to simplify the notation defined only for a special case we need).

\begin{de}
Let $y$ be a real number generator (ie. assume $R(y)$), $A$ is the set of $NE^*$-quotients of $u$ and $v$ with respect to $\ph\equiv y = 0\vv y\not = 0$ (then $\neg\ph$ is false ) if 
$$\forall x (\neg xv = u\ar \ph)\leftrightarrow(x \in A \vv \ph)$$
\end{de} 

\begin{lem}\label{NEQ}
Let $\ph$ be a formula as in the definition above, $A$ a species of natural numbers, $\alpha$ the corresponding choice sequence. Then RR-KS implies that there are global real number generators $u$ and $v$ such that $A$ is the set of $NE^*$-quotients of $u$ and $v$ with respect to $\ph$. 

\end{lem}
\proof
Let $\beta$ be the choice sequence given by RR-KS$(\ph,\alpha)$. 
\smallskip

We define two Cauchy sequences $u(n)$ and $v(n)$ as follows.
If $\forall l\leq n(\beta(l)=0)$, let $u(n)=v(n)=0$.
If $m=\mu n(\beta(n)>0)$, then for all $n\geq m$ let $u(n)=1/m$ and 
$v(n)=1/m\beta(m)$.
\smallskip

RR-KS$(\ph,\alpha)$ implies that $u$ and $v$ defined by these sequences have the required properties.
\qed

\bigskip

Let $L_{HAS}$ be the two sorted language of ${\bf HAS}$ with species variables but no sequence variables, $L = L_1\cup L_{HAS}$ and $L_2=L_{HAS}$ augmented with countably many constant symbols of species sort to be used in deductions.

\medskip

From now on let  $\ph\equiv y = 0\vv y\not = 0$ where the variable $y$ occurs only where it is indicated. We define a translation $\tau$ from the two sorted language $L_2$  into $L'$, the language of ordered rings augmented with countably many constant symbols - a pair $a_i$, $b_i$ corresponds to a constant symbol $A_i$ of $L_2$. We also fix a mapping from the set of variables of species sort into the set of pairs of variables of $L'$, mapping $X_i$ to $u_i$, $v_i$.

\medskip

Let $\Cal T'$ be the $L$-theory of universal closures of 
HAS $\cup$ RR-KS$(\ph,\alpha)$. 
\medskip

In an upcoming paper we prove the following.

\begin{prop}
RR-KS$(\ph,\alpha)$ holds in the model of \cite{sco1}. In particular $\Cal T'$ is consistent.\qed
\end{prop}

\noindent Since the schema R-KS is a special case of RR-KS, we can use the results of \cite{eszmjsl}. 
 In particular let the variables $x$, $y$, $u$, $v$, $w$ etc. range over reals,
let  $B(y)$ denote the L-formula $y = 0 \vee y \neq 0$ and let $\ph_\N(x,y,u,v) \equiv $
\[\neg x < 1 \wedge (\neg v = u \vee \neg xv = u
\rightarrow B(y)) \wedge \forall w [(\neg wv = u \rightarrow B(y))\rightarrow\]
\[ [(w < 1 \rightarrow B(y)) \wedge (w > 1 \rightarrow \exists w'(w \neq w' 
\vee \neg w'v = u + v \rightarrow B(y)))]].\]

\smallskip

\noindent Then the set of natural numbers is defined by an $L$-formula $\psi(x)\equiv \forall y\exists u\exists v\ph_\N(x,y,u,v)$, 
i.e. $\vdash\exists k\in\N(x=k)\leftrightarrow\psi(x)$ in two-sorted intuitionistic predicate calculus with equality, see Theorem 1. on page 1015  in \cite{eszmjsl}. 
We shall use the defined quantifiers
$$\en x\,\theta:\equiv  \er x(\psi(x)\wedge \theta)$$
$$\fn x\,\theta:\equiv  \fr x(\psi(x)\ar \theta)$$

\begin{de} For each $L_2$-formula $\theta$ we define its $L'$ translation $\tau(\theta)$ as follows.
\begin{itemize}
\item[-] If $\theta$ is first order atomic, $\tau(\theta):=\theta\vee\ph$;  in particular $\tau(\bot) := \ph$
\item[-] $\tau(n\in X_i) := \neg nu_i = v_i\ar\ph$
\item[-] $\tau(n\in A_i) := \neg na_i = b_i\ar\ph$
\item[-] $\tau(X_1 = X_2) := \tau(\forall x (x\in X_1\leftrightarrow x\in X_2))$
\item[-] $\tau(\psi_1\circ \psi_2) := \tau(\psi_1)\circ \tau(\psi_2)$ for $\circ = \wedge, \vee, \ar$
\item[-] $\tau(\neg\psi_1) := \tau(\psi_1)\ar \ph$ (special case of $\psi_1\ar \psi_2$)
\item[-] $\tau(\exists x\psi_1):=\en x(\tau(\psi_1))$
\item[-] $\tau(\forall x\psi_1):=\fn x(\tau(\psi_1))$
\item[-] $\tau(\exists X_i\psi_1:=\er u_i\er v_i \tau(\psi_1)$
\item[-] $\tau(\forall X_i\psi_1:=\fr u_i\fr v_i \tau(\psi_1)$

\end{itemize}
\end{de}

\medskip

\begin{lem}\label{addor}
For  an $L_2$-formula $\theta$,  \,
$\vdash \ph\ar\tau(\theta)$
\end{lem} 
\proof
 By formula induction.
\qed

\begin{lem}\label{deduction}
Let $\Sigma\cup \{\theta\}$ be a set of $L_2$ formulas (may containing free variables and constants), $\tau(\Sigma)$ the set of the translations. If $\Sigma\vdash\theta$ then $\tau(\Sigma)\vdash \tau(\theta)$.
\end{lem}

\proof By induction on the height of the natural deduction tree of $\Sigma\vdash\theta$. 

\begin{lem}\label{comp}
Let $\theta\equiv \exists X_i\forall n(n \in X_i\leftrightarrow A(x))$ be an instance of the comprehension schema CS where $A$ is an $L_{HAS}$-formula not containing  the species variable $X_i$ free. Then $\Cal T'\vdash \tau(\theta)$.
\end{lem}

\proof
$\tau(\theta)\equiv \er u_i\er v_i\fn x((\neg xv_i=u_i\ar\ph)\leftrightarrow \tau(A(x))$ 
By the comprehension schema CS applied to the formula  $\tau(A(x))$
 $$\Cal T'\vdash \exists X_i\forall x(x \in X_i \leftrightarrow \tau(A(x))$$

By Lemma \ref{NEQ}  $\Cal{T'}\vdash \forall X_i\er u_i\er v_i\fn x\bigl((\neg xv_i=u_i\ar\ph)\leftrightarrow x \in X_i\vee\ph\bigr)$.

If $\tau(A(x))$ holds, $x\in X_i$ and $\neg xv_i=u_i\ar\ph$ follows, so the right to left direction of $\tau(\theta)$ is true. If $\neg xv_i=u_i\ar\ph$, then either $x\in X_i$, so $\tau(\theta)$ holds, or $\ph$ holds. By Lemma \ref{addor} in both cases we have $ \tau(\theta)$.
\qed

\begin{lem}\label{choice}
Let $\theta\equiv \forall x\exists X_i B(X_i,x)\ar \exists X_i\forall xB((X_i)_x,x)$ be an instance of the choice schema AC-NF where $B$ is an $L_{HAS}$-formula. Then $\Cal T'\vdash \tau(\theta)$.
\end{lem}

\proof
$\tau(\theta)\equiv \fn x \,\er u_i\,\er v_i \,\tau(B)(u_i, v_i, x)\ar \er u_i\,\er v_i \, \fn x\,\tau(B)((u_i)_x, (v_i)_x, x)$. Two applications of AC-NF -- first, to $\er v_i \,\tau(B)(u_i, v_i, x)$ and then to $\tau(B)((u_i)_x, v_i, x)$ gives $\tau(\theta)$.
\qed

Since the translation of the axiom of extensionality $X_1 = X_2\leftrightarrow(\forall x (x\in X_1\leftrightarrow x\in X_2))$ is trivial, we have the following.

\begin{thm}
Let $\Cal{TH}$ be the set of universal closures of the formulas in $\tau \hbox{({\bf HAS})}$. For any $L_{\bf HAS}$- formula $\theta$, ${\bf HAS}$  $\vdash\theta$ implies $\Cal{TH}\vdash\tau(\theta)$ (so, since $\Cal{T'}\vdash \Cal{TH}$, $\Cal{T'}\vdash\tau(\theta)$) and ${\bf HAS}$  $\vdash\neg\theta$ implies $\Cal{TH}\not\vdash\tau(\theta)$.
\end{thm}
\proof
By Lemma \ref{deduction}  ${\bf HAS}$  $\vdash\theta \Rightarrow\tau \hbox{(${\bf HAS}$)}\vdash\tau(\theta)$. So if ${\bf HAS}$  $\vdash\neg\theta$, then   $\tau \hbox{({\bf HAS})}\vdash\tau(\theta)\ar(y=0\vv y\neq 0)$. From this, if $\Cal{TH}\vdash\tau(\theta)$, then $\Cal{TH}\vdash y=0\vee y\neq 0$. Then $\Cal{TH}\vdash\bot$ would follow, contradicting the fact that $\Cal{ TH}$ is consistent (since follows from a consistent theory $\Cal{T'}$). So in this case  $\Cal{TH}\not\vdash\tau(\theta)$.
\qed

\end{document}